\documentclass[11pt]{article}
\usepackage{amsmath}
\usepackage{amsthm}
\usepackage{amssymb}
\usepackage{geometry,mathtools}

\newtheorem{thm}{Theorem}

\newtheorem{lemma1}[thm]{Lemma}

\begin{document}
\date{}

\title{
Complete minors and average degree -- a short proof
}

\author{Noga Alon
\thanks{Department of Mathematics, Princeton University,
Princeton, NJ 08544, USA and
Schools of Mathematics and
Computer Science, Tel Aviv University, Tel Aviv 6997801,
Israel.
Email: {\tt nalon@math.princeton.edu}.
Research supported in part by
NSF grant DMS-1855464 and USA-Israel BSF grant 2018267.}
\and
Michael Krivelevich
\thanks{
School of Mathematical Sciences,
Tel Aviv University, Tel Aviv 6997801, Israel.
Email: {\tt krivelev@tauex.tau.ac.il}.  Research supported in part
by USA-Israel BSF grant 2018267.}
\and
Benny Sudakov
\thanks{
Department of Mathematics, ETH, Z\"urich, Switzerland.
Email: {\tt benjamin.sudakov@math.ethz.ch}.
Research supported in part by SNSF grant 200021\_196965.
}
}

\maketitle
\begin{abstract}
We provide a short and self-contained proof of the classical result
of Kostochka and of Thomason, ensuring that every graph of
average degree $d$ has a complete minor of order $\Omega(d/\sqrt{\log d})$.
\end{abstract}



Let $G=(V,E)$ be a graph with $|E|/|V|\ge d$. How large a complete minor are we guaranteed to find in $G$? This classical question, closely related to the famed Hadwiger's conjecture, has been thoroughly studied over the last half a century. It is quite easy to see the answer is at least logarithmic in $d$. Mader \cite{M68} proved it is of order at least  $d/\log d$. The right order of magnitude was established independently by Kostochka \cite{K82, K84} and by Thomason \cite{T84} to be $d/\sqrt{\log d}$, its tightness follows by considering random graphs. Finally, Thomason found in \cite{T01} the asymptotic value of this extremal function.

Here we provide a short and self-contained proof of the celebrated Kostochka--Thomason bound.

\begin{thm}
\label{th1}
Let $G=(V,E)$ be a graph with $|E|/|V|\ge d$, where $d$ is a sufficiently large integer. Then $G$ contains a minor of the complete graph on at least $\frac{d}{10\sqrt{\ln d}}$ vertices.
\end{thm}

The constant $1/10$ in the above statement is inferior to the best constant $3.13\ldots$ found by Thomason \cite{T01} (yet is better than the constants in \cite{K82, K84}); we did not make any serious attempt to optimize it in our arguments. The main  point here is to give a short proof of the tight $\Omega(d/\sqrt{\log d})$ bound for this classical extremal problem.

Throughout the proof we assume, whenever this is needed, that the
parameters $n$ and $d$ are sufficiently large. To simplify the presentation
we omit all floor and ceiling signs in several places. For a graph $G=(V,E)$, its minimum degree is denoted by $\delta(G)$, and for $v\in V$ we use $N_G(v)$ for the external neighborhood of $v$ in $G$.

We need the following lemma proven by simple probabilistic arguments. 
\begin{lemma1}\label{le1}
Let $H=(V,E)$ be a graph on at most $n$ vertices with $\delta(H)\ge n/6$. Let $t\le {n}/{\sqrt{\ln n}}$, and let $A_1,\ldots,A_t \subset V$ with  $|A_j|\le {n}{e^{-\sqrt{\ln n}/3}}$ for all $1\le j\le t$.
Then there is $B\subset V$ of size  $|B|\le 3.1\sqrt{\ln n}$ such that $B$ dominates all but at most ${n}{e^{-\sqrt{\ln n}/3}}$ vertices of $V$, $B\setminus A_j\ne\emptyset$ for all $j=1,\ldots,t$, and the induced subgraph $G[B]$ has at most six connected components.
\end{lemma1}
\begin{proof}
Set $s=3.1\sqrt{\ln n}$ and choose $s$ vertices of $V$ independently at random with repetitions. Let $B$ be the set of chosen vertices. Observe that for every vertex $v\in V$,
$$
Pr[N(v)\cap B=\emptyset]\le\left(1-\frac{d(v)}{n}\right)^s\le e^{-\frac{sd(v)}{n}}\le e^{-s/6}\,.
$$
Hence the expected number of vertices not dominated by $B$ is at
most $n e^{-s/6}<n e^{-3.1\sqrt{\ln n}/6}<n e^{-\sqrt{\ln n}/2}$, and  by Markov's inequality,
it is at most $n e^{-\sqrt{\ln n}/3}$ with probability exceeding $1/2$
(with room to spare). Also, since $|V| > \delta(H) \geq n/6$,
for every subset $A_j$,
$$
Pr[B\subseteq A_j]=
\left(\frac{|A_j|}{|V|}\right)^s <
\left(\frac{6|A_j|}{n}\right)^s
\le 6^s e^{-s\sqrt{\ln n}/3} = 6^{\Theta(\sqrt{\log n})}e^{-3.1\ln n/3}< \frac{1}{n}\,.
$$
Therefore the probability that $B\setminus A_j\ne\emptyset$ for all $j$
is at least $1-t/n\ge 1-1/\sqrt{\ln n}$.

We now argue about the number of connected components in $G[B]$. Writing $B=(v_1\ldots,v_s)$, for $1\le i\le s$ let  $x_i$ be the random variable counting the number of indices $1\le j\ne i\le s$ for which $v_j$ is a neighbor of $v_i$. Conditioning on $v_i$, we see that $x_i$ is distributed as a binomial random variable with parameters $s-1$ and  $d(v_i)/|V|>1/6$. Hence invoking the standard Chernoff-type bound on the lower tail of the binomial distribution, we derive that $Pr[x_i<s/7]\le e^{-\Theta(s)}$. Applying the union bound over all $1\le i\le s$, we conclude that with probability $1-o(1)$, we have $x_i\ge s/7$ for all $i$. Finally, observe that since $s\ll \sqrt{|V|}$, with probability $1-o(1)$ there are no repetitions in $B$, and hence $d(v_i,B)=x_i\ge s/7$ for all $1\le i\le s$. But then all connected components of $G[B]$ are of size exceeding $s/7$, and therefore $G[B]$ has at most six connected components.

Combining the above three estimates, the desired result follows.
\end{proof}

\noindent
{\bf Proof of Theorem \ref{th1}.}\,
 Let $G'=(V',E')$ be a minor of $G$ such that $|E'|\ge d|V'|$ and $|V'|+|E'|$ is minimal. If an edge $e$  of $G'$ is contained in $t$ triangles then contracting $e$ gives a minor of $G$ with one vertex and $t+1$ edges less. By the minimality of $G'$ we have $t+1>d$, implying $t\ge d$. Hence for every edge $e=(u,v)\in E(G')$, the vertex $u$ is connected by an edge of $G'$ to at least $d$ neighbors of $v$. The minimality of $G'$ also implies $|E'|=d|V'|$, hence $G'$ has a vertex $v$ of degree at most $2d$. Let $H$ be the subgraph of $G'$ induced by $N_{G'}(v)$. Then $H$ has at most $2d$ vertices and minimum degree at least $d$. Obviously a minor of $H$ is a minor of $G$ as well.

We now argue that $H$ contains a $d/3$-connected subgraph $H_1$ with $\delta(H_1)\ge 2d/3$.
If $H$ itself is $d/3$-connected this holds for $H_1=H$.
Otherwise there is a partition $V(H)=A\cup B\cup S$, where $A,B\ne\emptyset$, $|S|<d/3$, and $H$ has no edges between $A$ and $B$. Assume without loss of generality $|A|\le |B|$. Then $|A|\le d$, and since $\delta(H)\ge d$, every vertex $v\in A$ has at least $2d/3$ neighbors in $A$, implying that every pair of vertices of $A$ has at least $d/3$ common neighbors in $A$. Hence the induced subgraph $H_1:=H[A]$ is $d/3$-connected, has at most $2d$ vertices and satisfies $\delta(H_1)\ge 2d/3$.

Set $i=1$ and repeat the following iteration $d/10\sqrt{\ln d}$ times. Let $H_i=(V_i,E_i)\subseteq H_1$ be the current graph, and suppose $A_1,\ldots,A_{i-1}$ are subsets of $V_i$ of cardinalities $|A_j|\le {2d}{e^{-\sqrt{\ln(2d)}/3}}$ (representing the non-neighbors of the previously found branch sets $B_j$ in $V_i$). We assume (and justify it later) that $H_i$ is connected and has $\delta(H_i)> d/3$. Then the diameter of $H_i$ is at most $14$, as on any shortest path $P=(v_0,v_1,\ldots)$ in $H_i$ the vertices at
positions divisible by three have pairwise disjoint neighborhoods. Since
$|V(H_i)|/\delta(H_i)<6$, the number of such neighborhoods is at most $5$, and therefore any shortest path has at most $15$ vertices. Applying Lemma \ref{le1} with $H:=H_i$, $n:=2d$, $t:=i-1$, and $A_1,\ldots, A_{i-1}$ (for the initial step $i=1$ there are no $A_j$'s to plug into Lemma \ref{le1} --- which of course does not hinder its application) we get a subset $B_i$ of
cardinality $|B_i|\le 3.1\sqrt{\ln(2d)}$ as promised by the lemma. We now turn $B_i$ into a connected set by adding few vertices of $H_i$ if necessary. Recall that $H_i[B_i]$ has at most six connected components. Connecting one of them by shortest paths in $H_i$ to all others and recalling that $H_i$ has diameter at most 14, we conclude that by
appending to $B_i$ all the vertices of these paths
we make it connected by adding to it at most
$13\cdot 5=65$ vertices.
Altogether we obtain a connected subset $B_i$ of cardinality $|B_i|\le (3.1+o(1))\sqrt{\ln(2d)}$, dominating all but at most ${2d}{e^{-\sqrt{\ln(2d)}/3}}$ vertices of $V_i$ and having a vertex outside every $A_j$ (these properties are preserved under vertex addition when making $B_i$ into a connected subset) --- meaning connected to every previous $B_j$. We now update $V_{i+1}:=V_i-B_i$, $A_i:=V_{i+1}-N_{H_i}(B_i)$, and $A_j:=A_j\cap V_{i+1}$, $j=1,\ldots, i-1$, and finally increment $i:=i+1$, set $H_i:=H[V_i]$, and proceed to the next iteration. The total number of vertices deleted in all iterations satisfies:
$$
\big|\cup_{i} B_i\big|\le \frac{d}{10\sqrt{\ln d}}\cdot
(3.1+o(1))\sqrt{\ln(2d)}<\frac{d}{3}\,,
$$
and since we started with the $d/3$-connected graph $H_1$ with $\delta(H_1)\ge 2d/3$, we indeed have that at each iteration the graph $H_i$ is connected and has $\delta(H_i)> d/3$.

After having completed all $d/10\sqrt{\ln d}$ iterations, we get a family of $d/10\sqrt{\ln d}$ branch sets $B_i$, all connected, and with an edge of $H_1$ between every pair of branch sets. Hence they form a complete minor of order $d/10\sqrt{\ln d}$ as promised. \hfill$\Box$



\end{document}